\documentclass{article}
\usepackage{amsfonts}
\usepackage{amssymb}

\newtheorem{theorem}{Theorem}[section]

\newtheorem{corollary}{Corollary}

\newtheorem{definition}{Definition}[section]

\newtheorem{lemma}{Lemma}[section]

\newtheorem{proposition}{Proposition}[section]

\newenvironment{proof}[1][Proof]{\noindent\textbf{#1.} }{\ \rule{0.5em}{0.5em}}
\input{tcilatex}

\begin{document}

\title{AUTOMORPHISMS OF THE CATEGORY OF THE FREE NILPOTENT GROUPS OF THE
FIXED CLASS OF NILPOTENCY.}
\author{{\Large A.Tsurkov} \\
%EndAName
\textit{Department of Mathematics and Statistics,} \\
\textit{Bar Ilan University,} \\
\textit{Ramat Gan, 52900, Israel.\bigskip }\\
\textit{Jerusalem College of Technology, }\\
\textit{21 Havaad Haleumi, }\\
\textit{Jerusalem, 91160, Israel.}\\
tsurkov@jct.ac.il}
\maketitle

\begin{abstract}
This research was motivated by universal algebraic geometry. One of the
central questions of universal algebraic geometry is: when two algebras have
the same algebraic geometry? For answer of this question (see \cite{Pl},\cite%
{Ts}) we must consider the variety $\Theta $, to which our algebras belongs,
the category $\Theta ^{0}$ of all finitely generated free algebras of $%
\Theta $ and research how the group $\mathrm{Aut}\Theta ^{0}$ of all the
automorphisms of the category $\Theta ^{0}$ are different from the group $%
\mathrm{Inn}\Theta ^{0}$ of the all inner automorphisms of the category $%
\Theta ^{0}$. An automorphism $\Upsilon $ of the arbitrary category $%
\mathfrak{K}$ is called inner, if it is isomorphic as functor to the
identity automorphism of the category $\mathfrak{K}$, or, in details, for
every $A\in \mathrm{Ob}\mathfrak{K}$ there exists $s_{A}^{\Upsilon
}:A\rightarrow \Upsilon \left( A\right) $ isomorphism of these objects of
the category $\mathfrak{K}$ and for every $\alpha \in \mathrm{Mor}_{%
\mathfrak{K}}\left( A,B\right) $ the diagram%
\[
\begin{array}{ccc}
A & \overrightarrow{s_{A}^{\Upsilon }} & \Upsilon \left( A\right) \\ 
\downarrow \alpha &  & \Upsilon \left( \alpha \right) \downarrow \\ 
B & \underrightarrow{s_{B}^{\Upsilon }} & \Upsilon \left( B\right)%
\end{array}%
\]%
is commutative. In the case when $\Theta $ is a variety of all groups we
have the classical results which let us resolve this problem by an indirect
way. In \cite{DF} proved that for every free group $F_{i}$ the group $%
\mathrm{Aut}\left( \mathrm{Aut}F_{i}\right) $ coincides with the group $%
\mathrm{Inn}\left( \mathrm{Aut}F_{i}\right) $, from this result in \cite{Fo}
was concluded that $\mathrm{Aut}\left( \mathrm{End}F_{i}\right) =\mathrm{Inn}%
\left( \mathrm{End}F_{i}\right) $ and from this fact by theorem of reduction 
\cite{BPP} it can be concluded that $\mathrm{Aut}\Theta ^{0}=\mathrm{Inn}%
\Theta ^{0}$. In the case when $\Theta =\mathfrak{N}_{d}$ is a variety of
the all nilpotent class no more then $d$ groups we know by \cite{Ks}, that
if the number of generators $i$ of the free nilpotent class $d$ group $%
NF_{i}^{d}$ are bigger enough than $d$, then $\mathrm{Aut}\left( \mathrm{Aut}%
NF_{i}^{d}\right) =\mathrm{Inn}\left( \mathrm{Aut}NF_{i}^{d}\right) $. But
we have no description of $\mathrm{Aut}\left( \mathrm{End}NF_{i}^{d}\right) $
and so can not use the theorem of reduction. In this paper the method of
verbal operations is used. This method was established in \cite{PZ}. In \cite%
{PZ} by this method was very easily proved that $\mathrm{Aut}\Theta ^{0}=%
\mathrm{Inn}\Theta ^{0}$ when $\Theta $ is the variety of all groups and the
variety of all abelian groups. In this paper we will prove, by this method,
that $\mathrm{Aut}\mathfrak{N}_{d}^{0}=\mathrm{Inn}\mathfrak{N}_{d}^{0}$ for
every $d\geq 2$.
\end{abstract}

\section{Introduction and methodology.\label{Introduction}}

\setcounter{equation}{0}

This research was motivated by universal algebraic geometry. One of the
central questions of universal algebraic geometry is: when two algebras have
the same algebraic geometry? For answer of this question (see \cite{Pl},\cite%
{Ts}) we must consider the variety $\Theta $, to which our algebras belongs,
the category $\Theta ^{0}$ of all finitely generated free algebras of $%
\Theta $ and research how the group $\mathrm{Aut}\Theta ^{0}$ of all the
automorphisms of the category $\Theta ^{0}$ are different from the group $%
\mathrm{Inn}\Theta ^{0}$ of the all inner automorphisms of the category $%
\Theta ^{0}$. An automorphism $\Upsilon $ of the arbitrary category $%
\mathfrak{K}$ is called inner, if it is isomorphic as functor to the
identity automorphism of the category $\mathfrak{K}$, or, in details, for
every $A\in \mathrm{Ob}\mathfrak{K}$ there exists $s_{A}^{\Upsilon
}:A\rightarrow \Upsilon \left( A\right) $ isomorphism of these objects of
the category $\mathfrak{K}$ and for every $\alpha \in \mathrm{Mor}_{%
\mathfrak{K}}\left( A,B\right) $ the diagram%
\[
\begin{array}{ccc}
A & \overrightarrow{s_{A}^{\Upsilon }} & \Upsilon \left( A\right) \\ 
\downarrow \alpha &  & \Upsilon \left( \alpha \right) \downarrow \\ 
B & \underrightarrow{s_{B}^{\Upsilon }} & \Upsilon \left( B\right)%
\end{array}%
\]%
is commutative. In the case when $\Theta $ is a variety of all groups we
have the classical results which let us resolve this problem by an indirect
way. In \cite{DF} proved that for every free group $F_{i}$ the group $%
\mathrm{Aut}\left( \mathrm{Aut}F_{i}\right) $ coincides with the group $%
\mathrm{Inn}\left( \mathrm{Aut}F_{i}\right) $, from this result in \cite{Fo}
was concluded that $\mathrm{Aut}\left( \mathrm{End}F_{i}\right) =\mathrm{Inn}%
\left( \mathrm{End}F_{i}\right) $ and from this fact by theorem of reduction 
\cite{BPP} it can be concluded that $\mathrm{Aut}\Theta ^{0}=\mathrm{Inn}%
\Theta ^{0}$. In the case when $\Theta =\mathfrak{N}_{d}$ is a variety of
the all nilpotent class no more then $d$ groups we know by \cite{Ks}, that
if the number of generators $i$ of the free nilpotent class $d$ group $%
NF_{i}^{d}$ are bigger enough than $d$, then $\mathrm{Aut}\left( \mathrm{Aut}%
NF_{i}^{d}\right) =\mathrm{Inn}\left( \mathrm{Aut}NF_{i}^{d}\right) $. But
we have no description of $\mathrm{Aut}\left( \mathrm{End}NF_{i}^{d}\right) $
and so can not use the theorem of reduction. In this paper the method of
verbal operations is used. This method was established in \cite{PZ}. In \cite%
{PZ} by this method was very easily proved that $\mathrm{Aut}\Theta ^{0}=%
\mathrm{Inn}\Theta ^{0}$ when $\Theta $ is the variety of all groups and the
variety of all abelian groups. In this paper we will prove, by this method,
that $\mathrm{Aut}\mathfrak{N}_{d}^{0}=\mathrm{Inn}\mathfrak{N}_{d}^{0}$ for
every $d\geq 2$.

We start the explanation of the method of verbal operations in the general
situation. We consider the variety $\Theta $ of one-sorted algebras. The
signature of our algebras we denote $\Omega $.

For the construction of the category $\Theta ^{0}$ we must fix a countable
set of symbols $X_{0}=\left\{ x_{1},x_{2},\ldots ,x_{n},\ldots \right\} $.
As objects of the category $\Theta ^{0}$ we consider all free algebras $%
W\left( X\right) $ of the variety $\Theta $ generated by the finite subsets $%
X\subset X_{0}$. Morphisms of the category $\Theta ^{0}$ are homomorphisms
of these algebras.

\begin{definition}
\label{str_stab}\textbf{\hspace{-0.08in}. }\textit{We call the automorphism }%
$\Phi $\textit{\ of the category }$\Theta ^{0}$\textit{\ \textbf{strongly
stable} if it fulfills these three conditions:}

\begin{enumerate}
\item[A1)] $\Phi $ preserves all objects of $\Theta ^{0}$,

\item[A2)] there exists a system of bijections $\left\{ s_{B}^{\Phi
}:B\rightarrow B\mid B\in \mathrm{Ob}\Theta ^{0}\right\} $ such that $\Phi $
acts on the morphisms of $\Theta ^{0}$ by these bijections, i. e.,%
\begin{equation}
\Phi \left( \alpha \right) =s_{B}^{\Phi }\alpha \left( s_{A}^{\Phi }\right)
^{-1}  \label{1.1}
\end{equation}%
for every $\alpha :A\rightarrow B$ ($A,B\in \mathrm{Ob}\Theta ^{0}$);

\item[A3)] $s_{B}^{\Phi }\mid _{X}=id_{X}$ for every $B=W\left( X\right) \in 
\mathrm{Ob}\Theta ^{0}$.
\end{enumerate}
\end{definition}

The variety $\Theta $ is called an IBN variety if for every free algebras $%
W\left( X\right) ,W\left( Y\right) \in \Theta $ we have $W\left( X\right)
\cong W\left( Y\right) $ if and only if $\left\vert X\right\vert =\left\vert
Y\right\vert $. In the \cite[Theorem 2]{PZ} proved, that if $\Theta $ is an
IBN variety of one-sorted algebras, then every automorphism $\Psi \in 
\mathrm{Aut}\Theta ^{0}$ can be decomposed: $\Psi =\Upsilon \Phi $, where $%
\Upsilon ,\Phi \in \mathrm{Aut}\Theta ^{0}$, $\Upsilon $ is an inner
automorphism and $\Phi $ is a strongly stable one. So, if we want to know
the difference of the group $\mathrm{Aut}\Theta ^{0}$ from the group $%
\mathrm{Inn}\Theta ^{0}$, we must study the strongly stable automorphisms of
the category $\Theta ^{0}$.

Before the explanation of the notion of the verbal operation we will
introduce the short notation, which will be widely used in this paper. In
this notation $k$-tuple $\left( c_{1},\ldots ,c_{k}\right) \in C^{k}$ ($C$
is an arbitrary set) we denote by single letter $c$ and we will even allow
ourself to write $c\in C$ instead $c\in C^{k}$ and to write "homomorphism $%
\alpha :A\ni a\rightarrow b\in B$" instead "homomorphism $\alpha
:A\rightarrow B$, which transforms $a_{i}$ to the $b_{i}$, where $1\leq
i\leq k$".

For every word $w\left( x\right) \in W\left( X\right) $, where $X=\left\{
x_{1},\ldots ,x_{k}\right\} $ and every $C\in \Theta $ we can define a $k$%
-ary operation $w_{C}^{\ast }\left( c\right) =w\left( c\right) $ (in full
notation $c=\left( c_{1},\ldots ,c_{k}\right) \in C^{k}$, $x=\left(
x_{1},\ldots ,x_{k}\right) \in \left( W\left( X\right) \right) ^{k}$) or,
more formal, $w_{C}^{\ast }\left( c\right) =\gamma _{c}\left( w\left(
x\right) \right) $, where $\gamma _{c}$ is a well defined homomorphism $%
\gamma _{c}:W\left( X\right) \ni x\rightarrow c\in C$ (in full notation:
homomorphism $\gamma _{c}:W\left( X\right) \rightarrow C$, which transforms $%
x_{i}$ to the $c_{i}$ for $1\leq i\leq k$). This operation we call the 
\textbf{verbal operation} induced on the algebra $C$ by the word $w\left(
x\right) \in W\left( X\right) $. A little more detailed discussion about
definition of the verbal operation you can see in \cite[Section 2.1]{Ts}.

Now we will consider two kinds of substances, which will be also important
for our method:

\begin{enumerate}
\item systems of bijections $\left\{ s_{B}:B\rightarrow B\mid B\in \mathrm{Ob%
}\Theta ^{0}\right\} $ which fulfills these two conditions:

\begin{enumerate}
\item[B1)] for every homomorphism $\alpha :A\rightarrow B$ ($A,B\in \mathrm{%
Ob}\Theta ^{0}$) the mappings $s_{B}\alpha s_{A}^{-1}$ and $s_{B}^{-1}\alpha
s_{A}$ is also a homomorphism;

\item[B2)] $s_{B}\mid _{X}=id_{X}$ for every $B=W\left( X\right) \in \mathrm{%
Ob}\Theta ^{0}$.
\end{enumerate}

\item systems of words $\left\{ w_{\omega }\left( x\right) \in W\left(
X_{\omega }\right) \mid \omega \in \Omega \right\} $ which fulfills these
two conditions:

\begin{enumerate}
\item[Op1)] $X_{\omega }=\left\{ x_{1},\ldots ,x_{k}\right\} $, where $k$ is
an arity of $\omega $, for every $\omega \in \Omega $;

\item[Op2)] for every $B=W\left( X\right) \in \mathrm{Ob}\Theta ^{0}$ there
exists an isomorphism $\sigma _{B}:B\rightarrow B^{\ast }$ (algebra $B^{\ast
}$ has the same domain as the algebra $B$ and its operations $\omega
_{B}^{\ast }$ are induced by $w_{\omega }\left( x\right) $ for every $\omega
\in \Omega $) such as $\sigma _{B}\mid _{X}=id_{X}$.
\end{enumerate}
\end{enumerate}

The system of bijections $\left\{ s_{B}^{\Phi }:B\rightarrow B\mid B\in 
\mathrm{Ob}\Theta ^{0}\right\} $, described in A2) and A3) of the definition
of the strongly stable automorphism fulfills conditions B1) and B2) with $%
s_{B}=s_{B}^{\Phi }$.

We take a system of bijections $S=\left\{ s_{B}:B\rightarrow B\mid B\in 
\mathrm{Ob}\Theta ^{0}\right\} $ which fulfills conditions B1) and B2). If
arity of $\omega \in \Omega $ is $k$, we take $X_{\omega }=\left\{
x_{1},\ldots ,x_{k}\right\} \subset X_{0}$. $A_{\omega }=W\left( X_{\omega
}\right) $ - free algebra in $\Theta $. We have that $\omega \left( x\right)
\in A_{\omega }$ ($x=\left( x_{1},\ldots ,x_{k}\right) $) so there exists $%
w_{\omega }\left( x\right) \in A_{\omega }$ such that%
\begin{equation}
w_{\omega }\left( x\right) =s_{A_{\omega }}\left( \omega \left( x\right)
\right) .  \label{1.2}
\end{equation}%
The system of words $\left\{ w_{\omega }\left( x\right) \in A_{\omega }\mid
\omega \in \Omega \right\} $ we denote $\mathfrak{W}\left( S\right) $. This
system fulfills condition Op1) by our construction, condition Op2) with $%
\sigma _{B}=s_{B}$ ($B\in \mathrm{Ob}\Theta ^{0}$) by \cite[Theorem 3]{PZ}.
For every $C\in \Theta $ we denote $\omega _{C}^{\ast }$ the verbal
operation induced on the algebra $C$ by $w_{\omega }\left( x\right) $. $%
C^{\ast }$ will be the algebra, which has the same domain as the algebra $C$
and its operations are $\left\{ \omega _{C}^{\ast }\mid \omega \in \Omega
\right\} $. By \cite[Proposition 3.1]{Ts} one can conclude that for every $%
C\in \Theta $ the algebra $C^{\ast }$ belongs to $\Theta $.

Contrariwise, if we have a system of words $W=\left\{ w_{\omega }\left(
x\right) \mid \omega \in \Omega \right\} $, which fulfills conditions Op1)
and Op2), then the isomorphisms $\sigma _{B}:B\rightarrow B^{\ast }$ ($B\in 
\mathrm{Ob}\Theta ^{0}$) are bijections. For the system of bijections $%
\left\{ \sigma _{B}:B\rightarrow B\mid B\in \mathrm{Ob}\Theta ^{0}\right\} $%
, which we denote $\mathfrak{S}\left( W\right) $, B1) fulfills by \cite[1.8,
Lemma 8]{Gr} and \cite[Corollary from Proposition 3.2]{Ts}; B2) fulfills by
construction of these bijections.

If we have a system of bijections $S=\left\{ s_{B}:B\rightarrow B\mid B\in 
\mathrm{Ob}\Theta ^{0}\right\} $ which fulfills conditions B1) and B2) then
we can define an automorphism $\Phi \left( S\right) =\Phi $\ of the category 
$\Theta ^{0}$: $\Phi $ preserves all objects of the category $\Theta ^{0}$
and acts on its morphisms according to formula (\ref{1.1}) with $s_{B}^{\Phi
}=s_{B}$. Obviously $\Phi $ fulfills conditions A1) - A3) with $\left\{
s_{B}^{\Phi }:B\rightarrow B\mid B\in \mathrm{Ob}\Theta ^{0}\right\} =S$,
i.e. $\Phi $ is a strongly stable automorphism. Actually, the two different
systems of bijections $S_{1}=\left\{ s_{1,B}^{\Phi }:B\rightarrow B\mid B\in 
\mathrm{Ob}\Theta ^{0}\right\} $ and $S_{2}=\left\{ s_{2,B}^{\Phi
}:B\rightarrow B\mid B\in \mathrm{Ob}\Theta ^{0}\right\} $ can provide by
formula (\ref{1.1}) the same action on homomorphisms and, so, the same
strongly stable automorphism of the category $\Theta ^{0}$.

\begin{proposition}
\label{str_stab_obtained}\textbf{\hspace{-0.08in}. }\textit{Every strongly
stable automorphism }$\Phi $\textit{\ of the category }$\Theta ^{0}$\textit{%
\ can be obtained as }$\Phi \left( \mathfrak{S}\left( W\right) \right) $%
\textit{\ where }$W=\left\{ w_{\omega }\left( x\right) \in A_{\omega }\mid
\omega \in \Omega \right\} $\textit{\ is a system of words which fulfills
conditions Op1) and Op2).}
\end{proposition}

\begin{proof}
Let $\Phi $ be a strongly stable automorphism of the category $\Theta ^{0}$.
The system of bijections $S=\left\{ s_{B}^{\Phi }:B\rightarrow B\mid B\in 
\mathrm{Ob}\Theta ^{0}\right\} $ fulfills conditions B1) and B2). Lets
consider $\left\{ w_{\omega }\left( x\right) \in A_{\omega }\mid \omega \in
\Omega \right\} =\mathfrak{W}\left( S\right) $. By \cite[Proposition 3.3]{Ts}%
, $\mathfrak{W}\left( S\right) $ fulfills conditions Op1) and Op2). By \cite[%
Proposition 3.4]{Ts}, $\mathfrak{S}\left( \mathfrak{W}\left( S\right)
\right) =S$, so $\Phi \left( \mathfrak{S}\left( \mathfrak{W}\left( S\right)
\right) \right) =\Phi \left( S\right) $. $\Phi \left( S\right) $ and $\Phi $
both preserve all objects of $\Theta ^{0}$ and act on the morphisms of $%
\Theta ^{0}$ by the same system of bijections $S$ according to formula (\ref%
{1.1}). Therefore they coincide.
\end{proof}

So, if we describe the strongly stable automorphisms of the category $\Theta
^{0}$, we must concentrate on finding out the systems of words, which
fulfill conditions Op1) and Op2). However, in describing this, we must
remember, that different systems of words which fulfill conditions Op1) and
Op2) can provide us the same automorphism, because different systems of
bijections can provide us the same automorphism.

This method we apply to the variety $\mathfrak{N}_{d}$ of the all nilpotent
class no more then $d$ groups. The variety $\mathfrak{A}_{2}$ of the abelian
groups with the exponent no more then $2$ is a subvariety of $\mathfrak{N}%
_{d}$ for every $d\in 
%TCIMACRO{\U{2115} }%
%BeginExpansion
\mathbb{N}
%EndExpansion
$. So, by the second theorem of Fudzivara \cite[III.7.6]{Ku} the variety $%
\mathfrak{N}_{d}$ is an IBN variety.

\section{Applying the method.\label{applying}}

\setcounter{equation}{0}

The free $i$-generated group in $\mathfrak{N}_{d}$ is denoted as $NF_{i}^{d}$%
. The group signature is $\left\{ 1,-1,\cdot \right\} $, where $1$ is a $0$%
-ary operation of the taking of the unit, $-1$ is an unary operation of
taking an inverse element and $\cdot $ is a binary operation of the
multiplication. So, by our method, we must find out the systems of the words 
\begin{equation}
\left\{ w_{1},w_{-1},w_{\cdot }\right\}  \label{2.1}
\end{equation}%
such that $w_{1}\in NF_{0}^{d}$, $w_{-1}\in NF_{1}^{d}$, $w_{\cdot }\in
NF_{2}^{d}$ (condition Op1) ) and condition Op2) fulfills for all $%
NF_{i}^{d}\in \mathrm{Ob}\mathfrak{N}_{d}^{0}$ ($i\in 
%TCIMACRO{\U{2115} }%
%BeginExpansion
\mathbb{N}
%EndExpansion
$). But, as it will be clear above, we can consider only the word $w\left(
x,y\right) \in NF_{2}^{d}$ ($x,y$ are generators of $NF_{2}^{d}$) which
fulfills this condition:

\begin{enumerate}
\item[Op$^{d}$)] If $\left( NF_{2}^{d}\right) ^{\ast }$ is a set $NF_{2}^{d}$
with the binary verbal operation induced by $w\left( x,y\right) \in
NF_{2}^{d}$ (this operation we denote as "$\ast $"), then in $\left(
NF_{2}^{d}\right) ^{\ast }$ are fulfilled all group axioms and exists an
isomorphism $\sigma _{d}:NF_{2}^{d}\rightarrow \left( NF_{2}^{d}\right)
^{\ast }$, such that $\sigma _{d}\left( x\right) =x$, $\sigma _{d}\left(
y\right) =y$ ($x,y$ are the generators of $NF_{2}^{d}$).
\end{enumerate}

\begin{proposition}
\label{word_form_first}\textbf{\hspace{-0.08in}. }\textit{If in }$\left(
NF_{2}^{d}\right) ^{\ast }$\textit{\ are fulfilled all group axioms then }$%
w\left( x,y\right) =xyg_{2}\left( x,y\right) $\textit{\ (}$g_{2}\left(
x,y\right) \in \gamma _{2}\left( NF_{2}^{d}\right) $\textit{, }$\gamma
_{i}\left( G\right) $\textit{\ is a }$i$\textit{-th group in the lowest
central series of the group }$G$\textit{), }$1_{\ast }=1$\textit{, }$a^{\ast
k}=a^{k}$\textit{, for every }$a\in NF_{2}^{d}$\textit{, }$k\in 
%TCIMACRO{\U{2124} }%
%BeginExpansion
\mathbb{Z}
%EndExpansion
$\textit{\ and }$\gamma _{i}\left( \left( NF_{2}^{d}\right) ^{\ast }\right)
\subset \gamma _{i}\left( NF_{2}^{d}\right) $\textit{\ for every }$i\in 
%TCIMACRO{\U{2115} }%
%BeginExpansion
\mathbb{N}
%EndExpansion
$\textit{\ (}$1_{\ast }$\textit{\ and }$a^{\ast k}$\textit{\ is the unit of
the }$\left( NF_{2}^{d}\right) ^{\ast }$\textit{\ and degree of the element }%
$a$\textit{\ according to the new operation).}
\end{proposition}

\begin{proof}
Every word in $NF_{2}^{d}$ can be written as $w\left( x,y\right)
=x^{t}y^{s}g_{2}\left( x,y\right) $ ($t,s\in 
%TCIMACRO{\U{2124} }%
%BeginExpansion
\mathbb{Z}
%EndExpansion
$). So we assume that $x\ast y=x^{t}y^{s}g_{2}\left( x,y\right) $. Then $%
1\ast 1=1$, so $1=1_{\ast }$. Therefore $x=x\ast 1=x^{t}$, so $t=1$.
Analogously $s=1$. For every $a\in NF_{2}^{d}$ and every $k_{1},k_{2}\in 
%TCIMACRO{\U{2124} }%
%BeginExpansion
\mathbb{Z}
%EndExpansion
$ we have $a^{k_{1}}\ast a^{k_{2}}=a^{k_{1}+k_{2}}$, because $g_{2}\left(
a^{k_{1}},a^{k_{2}}\right) $ is calculated in the commutative group $%
\left\langle a\right\rangle $.

So, for every $a,b\in NF_{2}^{d}$ holds $\left( a,b\right) _{\ast
}=a^{-1}b^{-1}g_{2}\left( a^{-1},b^{-1}\right) \ast abg_{2}\left( a,b\right)
=a^{-1}b^{-1}abl_{2}=\left( a,b\right) l_{2}\in \gamma _{2}\left(
NF_{2}^{d}\right) $ ($l_{2}\in \gamma _{2}\left( NF_{2}^{d}\right) $, $%
\left( a,b\right) $ is a commutator of $a$ and $b$ and $\left( a,b\right)
_{\ast }$ is a commutator of $a$ and $b$ according to the new operation "$%
\ast $"). We assume that for $k<i$ it is proved that $\gamma _{k}\left(
\left( NF_{2}^{d}\right) ^{\ast }\right) \subset \gamma _{k}\left(
NF_{2}^{d}\right) $. If $l_{i-1}\in \gamma _{i-1}\left( \left(
NF_{2}^{d}\right) ^{\ast }\right) $, $b\in NF_{2}^{d}$, then we have 
\[
\left( l_{i-1},b\right) _{\ast }=l_{i-1}^{-1}b^{-1}g_{2}\left(
l_{i-1}^{-1},b^{-1}\right) \ast l_{i-1}bg_{2}\left( l_{i-1},b\right) = 
\]%
\[
=l_{i-1}^{-1}b^{-1}g_{2}\left( l_{i-1}^{-1},b^{-1}\right)
l_{i-1}bg_{2}\left( l_{i-1},b\right) g_{2}\left(
l_{i-1}^{-1}b^{-1}g_{2}\left( l_{i-1}^{-1},b^{-1}\right)
,l_{i-1}bg_{2}\left( l_{i-1},b\right) \right) . 
\]%
$g_{2}\left( l_{i-1}^{-1},b^{-1}\right) $, $g_{2}\left( l_{i-1},b\right) \in
\gamma _{i}\left( NF_{2}^{d}\right) $, because $l_{i-1}\in \gamma
_{i-1}\left( NF_{2}^{d}\right) $. So 
\[
\left( l_{i-1},b\right) _{\ast }\equiv l_{i-1}^{-1}b^{-1}l_{i-1}bg_{2}\left(
l_{i-1}^{-1}b^{-1},l_{i-1}b\right) \left( \mathrm{\func{mod}}\gamma
_{i}\left( NF_{2}^{d}\right) \right) . 
\]%
\[
\left( l_{i-1}^{-1}b^{-1},l_{i-1}b\right) \equiv \left(
l_{i-1}^{-1},b\right) \left( b^{-1},l_{i-1}\right) \equiv 1\left( \mathrm{%
\func{mod}}\gamma _{i}\left( NF_{2}^{d}\right) \right) . 
\]%
$g_{2}\left( x,y\right) $ has the weight $2$ or more, so $g_{2}\left(
l_{i-1}^{-1}b^{-1},l_{i-1}b\right) \in \gamma _{i}\left( NF_{2}^{d}\right) $
and $\left( l_{i-1},b\right) _{\ast }\in \gamma _{i}\left( NF_{2}^{d}\right) 
$.
\end{proof}

\begin{corollary}
\label{unit_inv}\textbf{\hspace{-0.08in}. }\textit{If the system of words (%
\ref{2.1})\ fulfills conditions Op1) and Op2) then }$w_{1}=1$\textit{, }$%
w_{-1}=x^{-1}$\textit{.}
\end{corollary}

By this Corollary we can concentrate on a research of the verbal binary
operations in the $NF_{2}^{d}$. We will find out the word $w\left(
x,y\right) \in NF_{2}^{d}$, which fulfill condition Op$^{d}$).

\begin{corollary}
\label{also_nilp}\textbf{\hspace{-0.08in}. }\textit{In the condition of the }%
Proposition \ref{word_form_first}\textit{\ the group }$\left(
NF_{2}^{d}\right) ^{\ast }$\textit{\ is also a nilpotent class }$d$\textit{\
group.}
\end{corollary}

\setcounter{corollary}{0}

By this Corollary, if in $\left( NF_{2}^{d}\right) ^{\ast }$\ are fulfilled
all group axioms then the homomorphism $\sigma _{d}:NF_{2}^{d}\ni \left( 
\begin{array}{c}
x \\ 
y%
\end{array}%
\right) \rightarrow \left( 
\begin{array}{c}
x \\ 
y%
\end{array}%
\right) \in \left( NF_{2}^{d}\right) ^{\ast }$ is well defined.

\begin{proposition}
\label{class2}\textbf{\hspace{-0.08in}. }\textit{The verbal operation
induced on }$NF_{2}^{2}$\textit{\ by the word }$w\left( x,y\right) =xy\left(
y,x\right) ^{m}$\textit{\ fulfills the group axioms for every }$m\in 
%TCIMACRO{\U{2124} }%
%BeginExpansion
\mathbb{Z}
%EndExpansion
$\textit{, but the homomorphism }$\sigma _{2}:NF_{2}^{2}\rightarrow \left(
NF_{2}^{2}\right) ^{\ast }$\textit{\ is an isomorphism if and only if }$m=0$%
\textit{\ or }$m=1$\textit{, i.e., }$w\left( x,y\right) =xy$\textit{\ or }$%
w\left( x,y\right) =yx$\textit{.}
\end{proposition}

\begin{proof}
Let $x\ast y=xy\left( y,x\right) ^{m}$. We have $1_{\ast }=1$, $a^{\ast
\left( -1\right) }=a^{-1}$.\textbf{\ }For every $a,b,c\in NF_{2}^{2}$ we have%
\[
\left( a\ast b\right) \ast c=ab\left( b,a\right) ^{m}\ast c= 
\]%
\[
=ab\left( b,a\right) ^{m}c\left( c,ab\left( b,a\right) \right)
^{m}=abc\left( b,a\right) ^{m}\left( c,a\right) ^{m}\left( c,b\right) ^{m}. 
\]%
\[
a\ast \left( b\ast c\right) =a\ast bc\left( c,b\right) ^{m}= 
\]%
\[
=abc\left( c,b\right) ^{m}\left( bc\left( c,b\right) ^{m},a\right)
^{m}=abc\left( c,b\right) ^{m}\left( b,a\right) ^{m}\left( c,a\right) ^{m}. 
\]%
So $\left( a\ast b\right) \ast c=a\ast \left( b\ast c\right) $.%
\[
\left( y,x\right) _{\ast }=y^{-1}\ast x^{-1}\ast y\ast x=y^{-1}x^{-1}\left(
x^{-1},y^{-1}\right) ^{m}\ast yx\left( x,y\right) ^{m}= 
\]%
\[
=\left( y,x\right) ^{1-2m}\left( yx,y^{-1}x^{-1}\right) ^{m}=\left(
y,x\right) ^{1-2m}\left( \left( y,x^{-1}\right) \left( x,y^{-1}\right)
\right) ^{m}=\left( y,x\right) ^{1-2m}. 
\]%
\[
\sigma _{2}\left( x^{\alpha _{1}}y^{\alpha _{2}}\left( y,x\right) ^{\alpha
_{3}}\right) =\sigma _{2}\left( x^{\alpha _{1}}\right) \ast \sigma
_{2}\left( y^{\alpha _{2}}\right) \ast \sigma _{2}\left( \left( y,x\right)
^{\alpha _{3}}\right) =x^{\alpha _{1}}\ast y^{\alpha _{2}}\ast \left( \left(
y,x\right) _{\ast }\right) ^{\alpha _{3}}= 
\]%
\[
=x^{\alpha _{1}}\ast y^{\alpha _{2}}\ast \left( \left( y,x\right) _{\ast
}\right) ^{\alpha _{3}}=x^{\alpha _{1}}\ast y^{\alpha _{2}}\ast \left(
y,x\right) ^{\left( 1-2m\right) \alpha _{3}}=x^{\alpha _{1}}y^{\alpha
_{2}}\left( y,x\right) ^{\alpha _{1}\alpha _{2}m+\left( 1-2m\right) \alpha
_{3}} 
\]%
($\alpha _{1},\alpha _{2},\alpha _{3}\in 
%TCIMACRO{\U{2124} }%
%BeginExpansion
\mathbb{Z}
%EndExpansion
$). If $1-2m\neq \pm 1$, i.e., $m\neq 0,1$, then $\left( y,x\right) \notin
\sigma _{2}\left( NF_{2}^{2}\right) $.
\end{proof}

\begin{proposition}
\label{gamma_i}\textbf{\hspace{-0.08in}. }\textit{If }$\left(
NF_{2}^{d}\right) ^{\ast }$\textit{\ as in the condition Op}$^{d}$\textit{),
then }$\gamma _{i}\left( NF_{2}^{d}\right) =\gamma _{i}\left( \left(
NF_{2}^{d}\right) ^{\ast }\right) $\textit{\ for every }$i\in 
%TCIMACRO{\U{2115} }%
%BeginExpansion
\mathbb{N}
%EndExpansion
$\textit{.}
\end{proposition}

\begin{proof}
By Corollary \ref{unit_inv} from Proposition \ref{word_form_first}, the
operations "$1$" and "$-1$" of the group signature in $NF_{2}^{d}$ and in $%
\left( NF_{2}^{d}\right) ^{\ast }$ coincide. By \cite[Proposition 3.2]{Ts},
original operations in $NF_{2}^{d}$ are the verbal operations induced on $%
NF_{2}^{d}$ by the words with respect to the operations "$1$", "$-1$" and "$%
\ast $". So we have a situation of Proposition \ref{word_form_first} in
which the group $NF_{2}^{d}$ and $NF_{2}^{d}$ changed they place. Hence, $%
\gamma _{i}\left( NF_{2}^{d}\right) \subset \gamma _{i}\left( \left(
NF_{2}^{d}\right) ^{\ast }\right) $ for every $i\in 
%TCIMACRO{\U{2115} }%
%BeginExpansion
\mathbb{N}
%EndExpansion
$.
\end{proof}

\begin{lemma}
\label{step}\textbf{\hspace{-0.08in}. }\textit{If the word }$w\left(
x,y\right) $\textit{\ fulfills conditions Op}$^{d}$\textit{), then the word }%
$\kappa \left( w\left( x,y\right) \right) \in NF_{2}^{d-1}$\textit{\
fulfills conditions Op}$^{d-1}$\textit{) (}$\kappa :NF_{2}^{d}\rightarrow
NF_{2}^{d-1}$\textit{\ is the natural epimorphism). }$\ $
\end{lemma}

\begin{proof}
$NF_{2}^{d-1}\in \mathfrak{N}_{d}$. So on $NF_{2}^{d-1}$ we can induce the
operations of the group signature by words $1\in NF_{0}^{d}$, $x^{-1}\in
NF_{1}^{d}$, $w\left( x,y\right) \in NF_{2}^{d}$. This operations we denote
by "$1$", "$-1$" and "$\circ $" correspondingly. $NF_{2}^{d-1}$ with this
operations we denote $\left( NF_{2}^{d-1}\right) ^{\circ }$. 
\begin{equation}
\begin{array}{ccccccccc}
1 & \rightarrow  & \gamma _{d}\left( NF_{2}^{d}\right)  & \hookrightarrow  & 
NF_{2}^{d} & \underrightarrow{\kappa } & NF_{2}^{d-1} & \rightarrow  & 1 \\ 
&  & \downarrow \sigma _{d} &  & \downarrow \sigma _{d} &  & \downarrow
\sigma _{d-1} &  &  \\ 
1 & \rightarrow  & \gamma _{d}\left( NF_{2}^{d}\right)  & \hookrightarrow  & 
\left( NF_{2}^{d}\right) ^{\ast } & \overrightarrow{\kappa } & \left(
NF_{2}^{d-1}\right) ^{\circ } & \rightarrow  & 1%
\end{array}
\label{2.2}
\end{equation}%
The operations in $\left( NF_{2}^{d}\right) ^{\ast }$ are induced by same
words; $\kappa :NF_{2}^{d}\rightarrow NF_{2}^{d-1}$ is a homomorphism, so,
by \cite[1.8, Lemma 8]{Gr}, $\kappa $ is also a homomorphism from $\left(
NF_{2}^{n}\right) ^{\ast }$ to $\left( NF_{2}^{n-1}\right) ^{\circ }$.
Therefore $\left( NF_{2}^{n-1}\right) ^{\circ }$ is an homomorphic image of
the group $\left( NF_{2}^{d}\right) ^{\ast }$, hence it is also a group. For
every $a,b\in NF_{2}^{d-1}$ we have $a\circ b=w\left( a,b\right) =\alpha
w\left( x,y\right) $, where $\alpha :NF_{2}^{d}\ni \left( 
\begin{array}{c}
x \\ 
y%
\end{array}%
\right) \rightarrow \left( 
\begin{array}{c}
a \\ 
b%
\end{array}%
\right) \in NF_{2}^{d-1}$. We can also consider a homomorphism $\widetilde{%
\alpha }:NF_{2}^{d-1}\ni \left( 
\begin{array}{c}
\kappa \left( x\right)  \\ 
\kappa \left( y\right) 
\end{array}%
\right) \rightarrow \left( 
\begin{array}{c}
a \\ 
b%
\end{array}%
\right) \in NF_{2}^{d-1}$. We have $\alpha =\widetilde{\alpha }\kappa $, so $%
a\circ b=\widetilde{\alpha }\left( \kappa w\left( x,y\right) \right) $.
Therefore the operation "$\circ $" is induced in $NF_{2}^{d-1}$ by the word $%
\kappa w=w\left( \kappa \left( x\right) ,\kappa \left( y\right) \right) \in
NF_{2}^{d-1}$ ($\kappa \left( x\right) $ and $\kappa \left( y\right) $ are
free generators of $NF_{2}^{d-1}$). This operation fulfills all group
axioms, hence, by Corollary \ref{also_nilp} from Proposition \ref%
{word_form_first}, $\left( NF_{2}^{d-1}\right) ^{\circ }$ is also a
nilpotent class $d-1$\ group. So $\sigma _{d-1}:NF_{2}^{d-1}\ni \left( 
\begin{array}{c}
\kappa \left( x\right)  \\ 
\kappa \left( y\right) 
\end{array}%
\right) \rightarrow \left( 
\begin{array}{c}
\kappa \left( x\right)  \\ 
\kappa \left( y\right) 
\end{array}%
\right) \in \left( NF_{2}^{d-1}\right) ^{\circ }$ is a well defined
homomorphism. Our goal is to prove that it is an isomorphism.

By Proposition \ref{gamma_i}, we have $\sigma _{d}\left( \gamma _{d}\left(
NF_{2}^{d}\right) \right) \subset \gamma _{d}\left( \left( NF_{2}^{d}\right)
^{\ast }\right) =\gamma _{d}\left( NF_{2}^{d}\right) $. By consideration of
the acting of $\sigma _{d}^{-1}$, we achieve $\gamma _{d}\left(
NF_{2}^{d}\right) \subset \sigma _{d}\left( \gamma _{d}\left(
NF_{2}^{d}\right) \right) $. In the diagram (\ref{2.2}) we have two exact
rows, because $\ker \kappa =\gamma _{d}\left( NF_{2}^{d}\right) $. The left
square of this diagram is commutative. The right square is commutative too,
because $\sigma _{d}$ and $\sigma _{d-1}$ preserve the generators of the
corresponding free nilpotent groups. $\sigma _{d}$ is an isomorphism, so $%
\sigma _{d-1}$ too.
\end{proof}

\begin{lemma}
\label{words2}\textbf{\hspace{-0.08in}. }\textit{For every }$d\geq 2$\textit{%
\ there are only two\ words: }$w=xy$\textit{\ and }$w=yx$\textit{\ from }$%
NF_{2}^{d}$\textit{\ which fulfills condition Op}$^{d}$\textit{).}
\end{lemma}

\begin{proof}
For $d=2$ it is proved in the Proposition \ref{class2}. We will assume that
it is proved for all natural numbers lesser then $d$.

Let $w\left( x,y\right) \in NF_{2}^{d}$ fulfills condition Op$^{d}$) and $%
\ast $ be an operation in $NF_{2}^{d}$ induced by $w\left( x,y\right) $. By
Lemma \ref{step} and assumption of induction $\kappa w\left( x,y\right)
=\kappa \left( x\right) \kappa \left( y\right) $ or $\kappa w\left(
x,y\right) =\kappa \left( y\right) \kappa \left( x\right) $.

In the first case we have $w\left( x,y\right) =xyr\left( x,y\right) $, where 
$r\left( x,y\right) \in \gamma _{d}\left( NF_{2}^{d}\right) $. For every $%
a,b,c\in NF_{2}^{n}$ we have $\left( a\ast b\right) \ast c=abr\left(
a,b\right) \ast c=abr\left( a,b\right) cr\left( abr\left( a,b\right)
,c\right) =abcr\left( a,b\right) r\left( ab,c\right) $ ($r\left( a,b\right)
=\alpha \left( r\left( x,y\right) \right) $ where $\alpha :NF_{2}^{d}\ni
\left( 
\begin{array}{c}
x \\ 
y%
\end{array}%
\right) \rightarrow \left( 
\begin{array}{c}
a \\ 
b%
\end{array}%
\right) \in NF_{2}^{d}$ for every $a,b\in NF_{2}^{d}$, so $r\left(
a,b\right) \in \gamma _{d}\left( NF_{2}^{d}\right) \subset Z\left(
NF_{2}^{d}\right) $). Also we have $a\ast \left( b\ast c\right) =a\ast
bcr\left( b,c\right) =abcr\left( b,c\right) r\left( a,bcr\left( b,c\right)
\right) =abcr\left( b,c\right) r\left( a,bc\right) $. $\ast $ fulfills the
axiom of associativity, so 
\begin{equation}
r\left( a,b\right) r\left( ab,c\right) =r\left( b,c\right) r\left(
a,bc\right) .  \label{2.3}
\end{equation}

We will prove that for $d\geq 3$ there isn't any nontrivial word $r\left(
x,y\right) $, which is generated by commutators of the weight $d$ and
fulfills the condition (\ref{2.3}). We will shift all our calculation from
the free nilpotent group $NF_{2}^{d}$ to the $NL_{2}^{d}$ - free nilpotent
class $n$ Lie algebra over $%
%TCIMACRO{\U{211a} }%
%BeginExpansion
\mathbb{Q}
%EndExpansion
$ with generators $x$ and $y$. By \cite[8.3.9]{Ba}, $\sqrt{NF_{2}^{d}}%
=\left( NL_{2}^{d}\right) ^{\circ }$, where $\sqrt{NF_{2}^{d}}$ is a Maltsev
completion of the group $NF_{2}^{d}$ and $\left( NL_{2}^{d}\right) ^{\circ }$
is a group which coincides with $NL_{2}^{d}$ as a set and has a
multiplication defined by Campbell-Hausdorff formula. By consideration of
the Campbell-Hausdorff formula we have that $r\left( a,b\right) =q\left(
a,b\right) $, where $q\left( x,y\right) $ is the word in $NL_{2}^{d}$, which
we achieve from $r\left( x,y\right) $ by replacement of the circular
brackets of the group commutators by Lie brackets and multiplication in $%
NF_{2}^{d}$ by addition in $NL_{2}^{d}$, $q\left( a,b\right) $ is a result
of substitution of the $a$ and $b$ instead $x$ and $y$ correspondingly ($%
a,b\in NF_{2}^{d}$) to the word $q\left( x,y\right) $. Also by consideration
of the Campbell-Hausdorff formula $ab\equiv a+b\left( \mathrm{\func{mod}}%
\left( NL_{2}^{d}\right) ^{2}\right) $ for every $a,b\in NF_{2}^{d}$ ($%
\left( NL_{2}^{d}\right) ^{i}$ is the $i$-th ideal of the lowest central
series of the algebra $NL_{2}^{d}$, $1\leq i\leq d$) and $ab=a+b$ if $a,b\in
\gamma _{d}\left( NF_{2}^{d}\right) $. Therefore the condition (\ref{2.3})
we can write as 
\begin{equation}
q\left( a,b\right) +q\left( a+b,c\right) =q\left( b,c\right) +q\left(
a,b+c\right) .  \label{2.4}
\end{equation}%
We shall substitute $a=\lambda x$, $b=x$, $c=y$ - ($\lambda \in 
%TCIMACRO{\U{2124} }%
%BeginExpansion
\mathbb{Z}
%EndExpansion
$) to (\ref{2.4}) and we will achieve:%
\begin{equation}
q\left( \left( \lambda +1\right) x,y\right) =q\left( x,y\right) +q\left(
\lambda x,x+y\right) .  \label{2.5}
\end{equation}%
$NL_{2}^{d}$ is a direct sum of its polyhomogeneous (homogeneous as
according $x$, as according $y$ separately) components (proof of this fact
is similar to the proof of \cite[2.2.5]{Ba}). Let $q\left( x,y\right)
=\sum\limits_{i=1}^{d-1}q_{i}$ $\left( x,y\right) $, where $q_{i}$ $\left(
x,y\right) $ is a polyhomogeneous component of $q\left( x,y\right) $
corresponding to the degree $i$ of $x$ and $d-i$ of $y$. By (\ref{2.5}) we
have 
\begin{equation}
\sum\limits_{i=1}^{d-1}\left( \lambda +1\right) ^{i}q_{i}\left( x,y\right)
=\sum\limits_{i=1}^{d-1}q_{i}\left( x,y\right)
+\sum\limits_{i=1}^{d-1}\lambda ^{i}q_{i}\left( x,x+y\right) .  \label{2.6}
\end{equation}%
When we develop $q_{i}\left( x,x+y\right) $ by additivity, we achieve $%
q_{i}\left( x,x+y\right) =q_{i}\left( x,y\right)
+\sum\limits_{j=i+1}^{d-1}m_{i,j}\left( x,y\right) $ where $m_{i,j}\left(
x,y\right) $ is a polyhomogeneous component of $q_{i}\left( x,x+y\right) $
with the degree $j$ of $x$ and $d-j$ of $y$. In particular $q_{1}\left(
x,x+y\right) =q_{1}\left( x,y\right) +\sum\limits_{j=2}^{d-1}m_{1,j}\left(
x,y\right) $. By comparison of the polyhomogeneous components of (\ref{2.6})
with the degree $2$ of $x$ and $d-2$ of $y$ we have%
\begin{equation}
\left( \left( \lambda +1\right) ^{2}-1-\lambda ^{2}\right) q_{2}\left(
x,y\right) =\lambda ^{2}m_{1,2}\left( x,y\right) .  \label{2.7}
\end{equation}

Let $\left\{ e_{1}\left( x,y\right) ,\ldots ,e_{k}\left( x,y\right) \right\} 
$ be a linear basis of the $d$-th ideal of the lowest central series of $%
NL_{2}^{d}$ and $q_{2}\left( x,y\right) =\sum\limits_{s=1}^{k}\rho
_{s}e_{s}\left( x,y\right) $, $m_{1,2}\left( x,y\right)
=\sum\limits_{s=1}^{k}\mu _{s}e_{s}\left( x,y\right) $ ($\rho _{s},\mu
_{s}\in 
%TCIMACRO{\U{2124} }%
%BeginExpansion
\mathbb{Z}
%EndExpansion
$). From (\ref{2.7}) we have%
\[
\left( \left( \lambda +1\right) ^{2}-1-\lambda ^{2}\right) \rho _{s}=\lambda
^{2}\mu _{s} 
\]%
($s\in \left\{ 1,\ldots ,k\right\} $). We take two values of $\lambda $ such
that we will achieve $\rho _{s}=0$, $\mu _{s}=0$. Hence $q_{2}\left(
x,y\right) =0$, $q_{2}\left( x,x+y\right) =0$, $m_{1,2}\left( x,y\right) =0$%
. Step by step, analogously we conclude, that $q_{i}\left( x,y\right) =0$, $%
q_{i}\left( x,x+y\right) =0$, $m_{1,i}\left( x,y\right) =0$ for $i\in
\left\{ 3,\ldots ,d-1\right\} $. Therefore $q\left( x,y\right) =q_{1}\left(
x,y\right) $. Now we shall substitute $a=x$, $b=y$, $c=\lambda y$ - ($%
\lambda \in 
%TCIMACRO{\U{2124} }%
%BeginExpansion
\mathbb{Z}
%EndExpansion
$) to (\ref{2.4}) and we will achieve: 
\[
q_{1}\left( x,y\right) +q_{1}\left( x+y,\lambda y\right) =q_{1}\left(
x,\left( \lambda +1\right) y\right) 
\]%
or%
\[
\left( \left( \lambda +1\right) ^{d-1}-\lambda ^{d-1}-1\right) q_{1}\left(
x,y\right) =0, 
\]%
so $q_{1}\left( x,y\right) =0$, $q\left( x,y\right) =0$, $r\left( x,y\right)
=1$ and $w\left( x,y\right) =xy$.

Analogously we consider the case $w\left( x,y\right) =yx$.
\end{proof}

\begin{theorem}
\label{main}\textbf{\hspace{-0.08in}. }\textit{All automorphisms of the
category of the nilpotent class }$d$\textit{\ free groups are inner.}
\end{theorem}

\begin{proof}
Let $W=\left\{ w_{1},w_{-1},w_{\cdot }\right\} $ be the system of words,
which fulfills conditions Op1) and Op2). By Lemma \ref{words2} there are no
more then two opportunities: $W=\left\{ 1,x^{-1},xy\right\} $ or may also be 
$W=\left\{ 1,x^{-1},yx\right\} $. It is easy to check that both these
systems actually fulfill conditions Op1) and Op2). By Proposition \ref%
{str_stab_obtained}, all strongly stable automorphisms of $\mathfrak{N}%
_{d}^{0}$ can by achieved as $\Phi \left( \mathfrak{S}\left( W\right)
\right) $. If $W$ fulfills conditions Op1) and Op2), then $\mathfrak{S}%
\left( W\right) $ fulfills conditions B1) and B2) (see Section \ref%
{Introduction}), so, by \cite[Lemma 3 and Theorem 2]{PZ}, $\Phi \left( 
\mathfrak{S}\left( W\right) \right) $ is an inner automorphism if and only
if for every $B\in \mathrm{Ob}\mathfrak{N}_{d}^{0}$ exists an isomorphism $%
c_{B}$ from $B$ to $B^{\ast }$, where $B^{\ast }$ is an algebra which has
the same domain as algebra $B$ and operations induced by the $\mathfrak{W}%
\left( \mathfrak{S}\left( W\right) \right) =W$ (see \cite[Proposition 3.4]%
{Ts}), such that $c_{B}\alpha c_{A}^{-1}=\alpha $ for every homomorphism $%
\alpha :A\rightarrow B$ ($A,B\in \mathrm{Ob}\Theta ^{0}$). If $W=\left\{
1,x^{-1},xy\right\} $, then $B^{\ast }=B$ and we can take $c_{B}=id_{B}$ for
every $B\in \mathrm{Ob}\mathfrak{N}_{d}^{0}$. And if $W=\left\{
1,x^{-1},yx\right\} $, then for every $B\in \mathrm{Ob}\mathfrak{N}_{d}^{0}$
we can take $c_{B}:B\rightarrow B^{\ast }$ such that $c_{B}\left( b\right)
=b^{-1}$ for every $b\in B$.
\end{proof}

\section{Acknowledgements.}

I dedicate this paper to the 80th birthday of Prof. B.Plotkin. He motivated
this research and was very heedful to it. Of course, I must appreciative to
Prof. B.Plotkin and Dr. G. Zhitomirski for theirs marvelous paper \cite{PZ}.
This research could not take place without theirs deep results. Prof. S.
Margolis also was very heedful to my research. Very useful discussion with
Prof. L. Rowen, Dr. R. Lipyanski, Dr. E. Plotkin and Dr. G. Zhitomirski,
also help me in the writing of this paper.

\end{document}